\documentclass{article}

\usepackage{amsmath,amssymb}
\usepackage{amsfonts}

\def\qed{$\square$}
\newtheorem{theorem}{Theorem}[section]

\newtheorem{example}[theorem]{Example}
\newtheorem{lemma}[theorem]{Lemma}
\newtheorem{proposition}[theorem]{Proposition}

\newtheorem{remark}[theorem]{Remark}

\begin{document}
\title{\bf Existence and smoothness of the density for spatially homogeneous SPDEs}
\author{
  \\
  {David Nualart}    \thanks{Supported by the NSF grant DMS-0604207  }         \\
  {\small \it Department of Mathematics }  \\[-0.15cm]
  {\small \it University of Kansas}          \\[-0.15cm]
  {\small \it Lawrence, Kansas, 66045, USA}  \\[-0.15cm]
  { \small {\tt nualart@math.ku.edu}}
  \\[-0.1cm]
\and
\\
  {Llu\'{\i}s Quer-Sardanyons \thanks{Supported by the grant MEC-FEDER Ref. MTM20006-06427 from the 
  Direcci\'on General de Investigaci\'on, Ministerio de Educaci\'on y Ciencia, Spain.}}              \\
  {\small\it Departament de Matem\`atiques} \\[-0.15cm]
  {\small\it  Universitat Aut\`onoma de Barcelona}          \\[-0.15cm]
  {\small\it 08193 Bellaterra (Barcelona),Spain}  \\[-0.15cm]
  {\small  {\tt quer@mat.uab.cat }} \\[-0.1cm]
  \\
}

\maketitle

\begin{abstract}
In this paper, we extend Walsh's stochastic integral with respect to a Gaussian noise, white in time and with some homogeneous spatial correlation, in order to be able to integrate some random measure-valued processes. This extension turns out to be equivalent to Dalang's one. Then we study existence and regularity of the density of the probability law for the real-valued mild solution to a general second order stochastic partial differential equation driven by such a noise. For this, we apply the techniques of the Malliavin calculus. Our results apply to the case of the stochastic heat equation in any space dimension and the stochastic wave equation in space dimension $d=1,2,3$. Moreover, for these particular examples, known results in the literature have been improved.

\end{abstract}

\vspace{2cm}


\vspace{0.3cm}


\newpage

\section{Introduction}

The purpose of this paper is to study the probability law of the real-valued solution of the following general class of stochastic partial differential equations:
\begin{equation}
L u(t,x) = \sigma(u(t,x)) \dot{W} (t,x) + b(u(t,x)), 
\label{1}
\end{equation}
$t \geq 0$, $x\in \mathbb{R}^{d}$, where $L$ denotes a second order differential operator, and we impose the
initial conditions
\begin{equation}
u(0,x) = \frac{\partial u}{\partial t}(0, x) = 0. \label{1a}
\end{equation}
The coefficients $\sigma$ and $b$ are some real-valued functions and $\dot{W} (t,x)$ is the formal notation for some Gaussian random perturbation 
defined on some probability space. We will assume that it is white in time and with a homogeneous spatial correlation given by a function $f$,  
and we denote by $(\mathcal{F}_t)_t$ the filtration generated by $W$ (see Section \ref{prel} for the precise definition of this noise).  \\

By definition, the solution to Equation (\ref{1}) is an adapted stochastic process $\{u(t,x), (t,x)\in [0,T]\times \mathbb{R}^{d}\}$
satisfying
\begin{align}
u(t,x)=& \int_0^t \int_{\mathbb{R}^{d}} \Gamma(t-s,x-y) \sigma(u(s,y)) W(ds,dy) \nonumber\\
&+ \int_0^t \int_{\mathbb{R}^{d}} b(u(t-s,x-y)) \Gamma(s,dy) ds,
\label{3}
\end{align}
where $\Gamma$ denotes the fundamental solution associated to $Lu=0$. As it is going to be clarified later on, for all $t\in [0,T]$, 
we suppose that $\Gamma(t)$ is a non-negative measure on $\mathbb{R}^d$ (see Hypothesis A).
The stochastic integral appearing in formula (\ref{3}) requires some care because the integrand
is a measure. For integrands that are real-valued functions, that is $(t,x)\mapsto \Gamma(t,x)\in \mathbb{R}$, the stochastic integral was defined by Walsh in \cite{walsh}. Then, in order to deal with SPDEs whose associated fundamental solution is a generalised function, Dalang \cite{Da} extended Walsh's stochastic integral and covered, for instance, the case of the wave equation in dimension three. \\

However, in the first part of this paper we will give, in a general setting, a definition of the stochastic integral with respect to $W$ using the techniques of the stochastic integration with respect to a cylindrical Brownian motion (see, for instance, \cite{dz}). As we will see, this integral will turn out to be equivalent to Dalang's extension (\cite{Da}) when the integrand is of the form $G:=\Gamma(t-\cdot,x-*) Z(\cdot,*)$, for certain stochastic processes $Z$. More precisely, we will show that random elements of this latter form may be integrated with respect to $W$ using a localising procedure: first we will assume that $Z$ has bounded trajectories and then we will identify $G$ as the weak limit of some sequence $(\Gamma Z_N)_N$, where any $Z_N$ have bounded paths, almost surely (see Lemma \ref{ggamma} and Proposition \ref{prop1}).\\

We should mention at this point that solutions to stochastic partial differential equations of the form (\ref{1}) have been largely studied during the last two decades. For instance, for the case of the stochastic heat and wave equations, we refer to \cite{walsh,carmona,DF,MS,Da,PZ,Pe}.\\

The second part of the paper is devoted to study the probability law of the random variable $u(t,x)$, for any fixed $(t,x)\in (0,T]\times \mathbb{R}^d$. This will be done using the techniques of the so-called Malliavin calculus. The aim is two-fold: \\

First, we prove that $u(t,x)$ has an absolutely continuous law with respect to Lebesgue measure on $\mathbb{R}$, provided that, among other assumptions, the differential operator $L$ and the spatial correlation $f$ are related as follows: 
\begin{equation}
\int_0^T \int_{\mathbb{R}^{d}} |\mathcal{F} \Gamma(t)(\xi)|^2 \mu(d\xi) dt< \infty,
\label{5bis}
\end{equation}
where $\mu$ is a non-negative tempered measure such that its Fourier transform is $f$
(see Theorem \ref{densitat} for the precise statement). This result provides a generalisation of Theorem 3 in \cite{QS}, where the authors deal with the three-dimensional stochastic wave equation and a slightly stronger condition than (\ref{5bis}) is assumed. In order to prove that the law of $u(t,x)$ has a density, we apply Bouleau-Hirsch's criterion. Indeed, to prove the Malliavin regularity of the solution, we take advantage of the results in 
\cite{QS} and we fully identify the initial condition of the stochastic equation satisfied by the Malliavin derivative of $u(t,x)$. This is an important point in order to show that the Malliavin matrix is invertible almost surely. Eventually, we point out that (\ref{5bis}) is also a sufficient condition to have the stochastic integral in (\ref{3}) well defined. \\

Secondly, we prove that the law of $u(t,x)$ has an infinitely differentiable density with respect to Lebesgue measure on $\mathbb{R}$ (see 
Theorem \ref{regularitatdensitat}). To obtain this result, we show that $u(t,x)$ is infinitely differentiable in the Malliavin sense and that the inverse of the Malliavin matrix has moments of all order. For the latter to be achieved, we need to impose a lower bound of the integral in (\ref{5bis})   in terms of a certain power of $T$ (see (\ref{lowbound})).
We should mention that 
Theorem \ref{regularitatdensitat} in Section \ref{regularitat} provides an improvement of Theorem 3 in \cite{qs2} for the case of the three-dimensional stochastic wave equation, since in this latter reference the integrability conditions concerning the Fourier transform of $\Gamma$ are much more involved (see also \cite{marta})). Moreover, it is worth mentioning that Theorem \ref{regularitatdensitat} also generalises known results on existence and smoothness of the density for the case of the stochastic heat and wave equation with dimensions $d\geq 1$ and $d=1,2$, respectively (see 
\cite{carmona,MS,mms,marta}).\\

The paper is organised as follows. In the next Section \ref{prel}, we present some preliminaries concerning the random perturbation in Equation (\ref{1}) as well as the main hypothesis on the fundamental solution $\Gamma$ and the space correlation $f$. We extend Walsh's stochastic integral in order to cover the case of measure-valued integrands in Section \ref{integral}; we also caracterise some measure-valued random elements that can be integrated with respect to $W$ and we sketch the construction of the integral in a Hilbert-valued setting. In Section \ref{ex}, a theorem on existence and uniqueness of solution for Equation (\ref{1}) is stated; we also deal with some particular examples of differential operators $L$, namely the heat and wave equations.   Sections \ref{existencia} and \ref{regularitat} are devoted, respectively, to the existence and regularity of the density for the probability law of the solution to (\ref{1}).  At the very beginning of Section \ref{existencia}, we introduce the main tools of the Malliavin calculus needed along the paper (we refer to \cite{nualart} for a complete account on the topic).\\

Throughout the paper we use the notation $C$ for any positive real constant, independently of its value.

\section{Preliminaries}
\label{prel}

Recall that we are interested in the following general class of stochastic partial differential equations:
$$
L u(t,x) = \sigma(u(t,x)) \dot{W} (t,x) + b(u(t,x)), 
$$
with $t \geq 0$, $x\in \mathbb{R}^{d}$, $L$ denotes a second order differential operator, and we consider vanishing initial conditions (\ref{1a}).
The Gaussian random perturbation is described as follows: $W$ is a zero mean Gaussian family of random variables  
$\{W(\varphi), \varphi \in \mathcal{C}_0^\infty (\mathbb{R}^{d+1})\}$,
defined in a complete probability space $(\Omega, \mathcal{F},P)$,  with covariance
\begin{equation}
E(W(\varphi) W(\psi)) = \int_0^\infty \int_{\mathbb{R}^{d}} \int_{\mathbb{R}^{d}} \varphi(t,x) f(x-y) \psi(t,y) dx dy dt,
\label{2}
\end{equation}
where $f$ is a non-negative continuous function of $\mathbb{R}^{d}\setminus \{0\}$ such that it is the Fourier
transform of a non-negative definite tempered measure $\mu$ on $\mathbb{R}^{d}$.
That is,
$$f(x)=\int_{\mathbb{R}^{d}} \exp(-2\pi i\; x\cdot \xi) \mu(d\xi)$$
and there is an integer $m\geq 1$ such that
$$\int_{\mathbb{R}^{d}} (1+|\xi|^2)^{-m}\mu(d\xi) <\infty.$$
Then, the covariance (\ref{2}) can also be written, using Fourier transform, as
$$E(W(\varphi) W(\psi)) = \int_0^\infty \int_{\mathbb{R}^{d}} \mathcal{F} \varphi(t)(\xi) \overline{\mathcal{F} \psi(t)(\xi)} \mu(d\xi) dt.$$
The main assumption on the differential operator $L$ may be summarised as follows:\\

\noindent {\bf Hypothesis A}.
The fundamental solution to $Lu = 0$, denoted by $\Gamma$, is a non-negative measure of
the form $\Gamma(t, dx)dt$ such that for all $T>0$
\begin{equation}
\sup_{0\leq t\leq T} \Gamma(t,\mathbb{R}^{d})\leq C_T<\infty
\label{e2}
\end{equation}
and 
\begin{equation}
\int_0^T \int_{\mathbb{R}^{d}} |\mathcal{F} \Gamma(t)(\xi)|^2 \mu(d\xi) dt< \infty.
\label{5}
\end{equation}

\vspace{0.5cm}

The completion of the Schwartz space $\mathcal{S} (\mathbb{R}^{d})$ of rapidly decreasing $\mathcal{C}^\infty$ functions, endowed
with the inner product
$$\langle \varphi,\psi\rangle_\mathcal{H}  
=\int_{\mathbb{R}^{d}} \int_{\mathbb{R}^{d}} \varphi(x) f(x-y) \psi(y) dx dy=\int_{\mathbb{R}^{d}} \mathcal{F} \varphi(\xi) \overline{\mathcal{F} \psi(\xi)} \mu(d\xi),$$
$\varphi, \psi\in \mathcal{S} (\mathbb{R}^{d})$, is denoted by $\mathcal{H}$. Notice that $\mathcal{H}$ may contain distributions. Set $\mathcal{H}_T = L^2([0,T];\mathcal{H})$.


\section{Stochastic integrals}
\label{integral}

In this section Walsh's stochastic integral with respect to martingale measures will be extended to more general integrands, namely the class of square integrable $\mathcal{H}$-valued predictable processes. The extension will be performed in the infinite dimensional setting described by Da Prato and Zabczyk in \cite{dz}. 
Then, we will give non-trivial examples of integrands, which will be some measure-valued random elements. We will briefly recall the extension of the stochastic integral in a Hilbert-valued setting. This will be needed to give a rigorous meaning to the stochastic evolution equations satisfied by the Malliavin derivatives of the solution of (\ref{3}).\\

\noindent {\bf Extension of Walsh's stochastic integral}\\

Fix a time interval $[0,T]$. The Gaussian family $\{W(\varphi ),\varphi \in
\mathcal{C}_{0}^{\infty }([0,T]\times \mathbb{R}^{d})\}$ can be extended to the
completion $\mathcal{H}_T=L^{2}([0,T];\mathcal{H})$ of the space $\mathcal{C}_{0}^{\infty
}([0,T]\times \mathbb{R}^{d})$ under the scalar product%
\[
\left\langle \varphi ,\psi \right\rangle =\int_{0}^{T}\int_{\mathbb{R}^{d}}%
\mathcal{F}\varphi (t)(\xi )\overline{\mathcal{F}\psi (t)(\xi )}\mu (d\xi
)dt. 
\]%
We will also denote by $W(g)$ the Gaussian random variable associated with
an element $g\in L^{2}([0,T];\mathcal{H})$.

Set $W_{t}(h)=W(1_{[0,t]}h)$ for any $t\geq 0$ and $h\in \mathcal{H}$. \
Then, $\{W_{t},t\in \lbrack 0,T]\}$ is a cylindrical Wiener process in the
Hilbert space $\mathcal{H}$. That is, for any $h\in \mathcal{H}$, $%
\{W_{t}(h),t\in \lbrack 0,T]\}$ is a Brownian motion with variance $\left\|
h\right\| _{\mathcal{H}}^{2}$, and%
\[
E(W_{t}(h)W_{s}(g)=\left( s\wedge t\right) \left\langle h,g\right\rangle _{%
\mathcal{H}}. 
\]%
Let $\mathcal{F}_t$ be the $\sigma$-field generated by the random variables
$\{W_s(h), h\in \mathcal{H}, 0\le s\le t\}$ and the $P$-null sets. We define the predictable
$\sigma$-field as the $\sigma$-field in $\Omega\times [0,T]$ generated by
the sets $\{ (s,t]\times A, 0\le s<t\le T, A\in \mathcal{F}_s\}$.

Then (see, for instance, \cite{dz}), we can define the
stochastic integral of $\mathcal{H}$-valued square integrable predictable
processes. For any predictable process $g\in $ $L^{2}(\Omega \times \lbrack
0,T];\mathcal{H)}$ we denote its integral with respect to the cylindrical
Wiener process $W$ by%
\begin{equation}
\int_{0}^{T}\int_{\mathbb{R}^{d}}gdW=g\cdot W,  \label{e1}
\end{equation}
and we have the isometry property%
\[
E\left( \left| g\cdot W\right| ^{2}\right) =E\left( \int_{0}^{T}\left\|
g_{t}\right\| _{\mathcal{H}}^{2}dt\right) . 
\]

\begin{remark}
Under the standing assumptions, using an approximation procedure by means of test functions,
one proves that the space $\mathcal{H}$ contains the indicator functions of bounded Borel sets (for details see \cite{DF} or 
\cite{quer}, p. 13). Then, $M_{t}(A):=W(1_{[0,t]}1_{A})$ defines a martingale measure associated to the noise $%
W $ in the sense of Walsh (see \cite{walsh} and \cite{Da}) and the stochastic integral  (\ref{e1})
coincides with the integral defined in the work of Dalang \cite{Da}.
\end{remark}

\noindent {\bf Example of integrands}\\

We aim now to provide useful examples of  random distributions which belong to the space $L^2(\Omega\times [0,T]; \mathcal{H})$. 
Before stating the result, we consider the following lemma:

\begin{lemma}\label{ggamma}
Assume that $\Gamma$ satisfies Hypothesis A. Let $g$ be a bounded Borel function on $[0,T]\times \mathbb{R}^{d}$. Then $%
g\Gamma \in \mathcal{H}_T$, and%
\[
\left\| g\Gamma \right\|_{\mathcal{H}_T}^2 \leq \left\| g\right\|^2
_{\infty }\int_{0}^{T}\int_{\mathbb{R}^{d}}\left| \mathcal{F}\Gamma (t)(\xi
)\right| ^{2}\mu (d\xi )dt.
\]
\end{lemma}

\noindent {\it Proof.}
We can decompose $g$ into the difference $g^{+}-g^{-}$ of two nonnegative
bounded Borel functions. Thus, without any loss of generality we can assume
that $g$ is nonnegative. Moreover, we observe that $g\Gamma $ also satisfies \ conditions (%
\ref{e2}) and (\ref{5}). Indeed, to prove the latter condition, we consider  
an approximation of the identity $(\psi _{n})_n$ defined as follows: let $\psi\in \mathcal{C}_0^{\infty}(\mathbb{R}^d)$ such that $\psi\geq 0$, 
the support of $\psi$ is contained in the unit ball of $\mathbb{R}^d$ and $\int_{\mathbb{R}^d} \psi(x)dx=1$. 
Set $\Gamma
_{n}(t)=\psi _{n}\ast \Gamma (t)$ and $J_{n}(t)=\psi _{n}\ast (g\Gamma (t))$. 
Then, for all $t\in \lbrack 0,T]$,  $\Gamma _{n}(t)$ and  $J_n(t)$ belong to $\mathcal{S}(\mathbb{R}^{d})\subset \mathcal{H}$, and \ $\left| \mathcal{F}%
\Gamma _{n}(t)\right| \leq \left| \mathcal{F}\Gamma (t)\right|$. Besides, since $\Gamma$ is non-negative, we have that 
$J_n(t)(\xi)\leq \|g\|_\infty \Gamma_n(t)(\xi)$, for any $\xi\in \mathbb{R}^d$.  Thus, by Fatou's lemma
\begin{align*}
& \int_0^T \int_{\mathbb{R}^{d}} |\mathcal{F} (g\Gamma(t))(\xi)|^2 \mu(d\xi) dt \\
&\quad \leq \liminf_{n\rightarrow \infty} \int_0^T  \int_{\mathbb{R}^{d}}\int_{\mathbb{R}^{d}} J_n(t,x) f(x-y) J_n(t,y) dx dy dt \\
&\quad \leq \|g\|_{\infty}^2 \liminf_{n\rightarrow \infty} \int_0^T \int_{\mathbb{R}^{d}} |\mathcal{F} \Gamma_n(t)(\xi)|^2 \mu(d\xi) dt\\
& \quad \leq \|g\|_{\infty}^2 \int_0^T \int_{\mathbb{R}^{d}} |\mathcal{F} \Gamma(t)(\xi)|^2 \mu(d\xi) dt < \infty. 
\end{align*}
The fact that $g\Gamma$ satisfies conditions (\ref{e2}) and (\ref{5}) let us reduce the proof to the case where $g=1$. 

We consider the regularisation $(\Gamma_n)_n$ of $\Gamma$ defined above. Condition (\ref{5}) implies that $\Gamma _{n}(t)$ belongs to 
$\mathcal{H}_T$ and it has a uniformly bounded norm, so it converges weakly to some
element $h\in \mathcal{H}_T$. We claim that $h=\Gamma $, and this is a consequence of the
fact that, owing to the definition of $\Gamma_n(t)$, for any $\varphi \in \mathcal{S}(\mathbb{R}^{d})$, and for any $%
0\leq s<t\leq T$ we have%
\[
\int_{s}^{t}\left\langle h,\varphi \right\rangle _{\mathcal{H}%
}dr=\lim_{n\rightarrow \infty }\int_{s}^{t}\left\langle \Gamma
_{n}(r),\varphi \right\rangle _{\mathcal{H}}dr=\int_{s}^{t}\left\langle
\Gamma (r),\varphi \right\rangle _{\mathcal{H}}dr.
\]
More precisely, it holds that
$$
\int_{s}^{t}\left\langle \Gamma
_{n}(r),\varphi \right\rangle _{\mathcal{H}}dr  =  \int_{s}^{t} \int_{\mathbb{R}^d} 
\Gamma(r,dz) \left( \int_{\mathbb{R}^d} \psi_n(x-z) F(x) dx \right)dr,$$
where $F(x):=\int_{\mathbb{R}^d} f(x-y)\varphi(y) dy$. Observe that the hypothesis on $f$ and $\varphi$ imply that $F$ is continuous and 
$\lim_{x\rightarrow \infty} F(x)=0$. Hence, it turns out that we can apply the Bounded Convergence Theorem, so that we end up with
$$
\lim_{n\rightarrow \infty } \int_{s}^{t}\left\langle \Gamma
_{n}(r),\varphi \right\rangle _{\mathcal{H}}dr   =  
\int_{s}^{t} \int_{\mathbb{R}^d} \int_{\mathbb{R}^d} \Gamma(r,dz)  f(z-y) \varphi(y) dy dr,$$
and this let us identify $h$ with $\Gamma$ in $\mathcal{H}_T$. \hfill
\qed
 
\medskip
This lemma allows us to prove the following result, which give examples of random distributions that can be integrated with respect to $W$.

\begin{proposition}\label{prop1} 
Assume that $\Gamma$ satisfies Hypothesis A.
Let $Z=\{Z(t,x),(t,x)\in \lbrack 0,T]\times \mathbb{R}^{d}\}$ be a
predictable process such that%
\[
C_{Z}:=\sup_{(t,x)\in \lbrack 0,T]\times \mathbb{R}^{d}}E(|Z(t,x)|^{2})<%
\infty .
\]%
Then, the random element $G=G(t,dx)=Z(t,x)\Gamma (t,dx)$ is a predictable
process in the space $L^{2}(\Omega \times \lbrack 0,T];\mathcal{H)}$.
\end{proposition}

\noindent {\it Proof.}
For any $N\geq 1$ define%
\[
Z_{N}(t,x)=Z(t,x){\bf 1}_{\{|Z(t,x)|\leq N\}}.
\]%
Clearly, $Z_{N}$ is a predictable process with bounded trajectories. Thus,
by the previous lemma, $G_{N}(t,x):=Z_{N}(t,x)\Gamma (t,dx)$ is a
predictable process in $L^{2}(\Omega \times \lbrack 0,T];\mathcal{H)}$.
Let $\left( J_{N,n}(t) \right)_n \subset \mathcal{S}(\mathbb{R}^d) $ be the regularisation of $G_N(t)$ by means of an approximation of the identity $(\psi_n)_n$, as it has been defined in the proof of Lemma \ref{ggamma}. This let us prove that the norm of $G_N$ in $L^2(\Omega\times [0,T];\mathcal{H})$
is uniformly bounded because%
\begin{align*}
& E\left( \left\| G_{N}\right\| _{L^{2}([0,T];\mathcal{H})}^{2}\right) 
= E\left( \lim_{n\rightarrow \infty} \|J_{N,n}\|_{L^{2}([0,T];\mathcal{H})}^{2}\right)\\
&\quad \leq \liminf_{n\rightarrow \infty} E\left( \int_0^T \int_{\mathbb{R}^{d}}\int_{\mathbb{R}^{d}} J_{N,n}(t,x) f(x-y)J_{N,n}(t,y) dx dy dt\right)\\
&\quad \leq C_Z \liminf_{n\rightarrow \infty}  E\left( \int_{0}^{T}\int_{\mathbb{R}^{d}}\int_{\mathbb{R}^{d}}\Gamma_n
(t,x) f(x-y) \Gamma_n(t,y) dx dy dt \right)\\
&\quad \leq C_{Z}\int_{0}^{T}\int_{\mathbb{R}^{d}}\left| \mathcal{F}\Gamma
(t)(\xi )\right| ^{2}\mu (d\xi )dt < \infty.
\end{align*}%
Recall that $\Gamma_n(t)=\psi_n * \Gamma(t)$. 
Therefore, $G_{N}$ converges weakly to a predictable process $\widetilde{G}$
in $L^{2}(\Omega \times \lbrack 0,T];\mathcal{H)}$. We claim that $%
\widetilde{G}=G$. In fact, for any $\varphi \in \mathcal{S}(\mathbb{R}^{d})$%
,  for any $0\leq s<t\leq T$ and for any $B\in \mathcal{F}_{s}$, we can argue similarly as in the very last part of the proof of Lemma \ref{ggamma} and obtain 
\begin{align*}
& E\left( {\bf 1}_{B}\int_{s}^{t}\left\langle \widetilde{G}(r),\varphi \right\rangle
_{\mathcal{H}}dr\right) \\
&\quad = \lim_{N\rightarrow \infty }E\left(
{\bf 1}_{B}\int_{s}^{t}\left\langle G_{N}(r),\varphi \right\rangle _{\mathcal{H}%
}dr\right)  \\
&\quad =\lim_{N\rightarrow \infty }E\left( {\bf 1}_{B}\int_{s}^{t}\left\langle
Z(r,x){\bf 1}_{\{|Z(r,x)|\leq N\}}\Gamma (r),\varphi \right\rangle _{\mathcal{H}%
}dr\right)  \\
&\quad =\lim_{N\rightarrow \infty }E\left( {\bf 1}_B \int_{s}^{t}\int_{\mathbb{R}^{d}}\int_{%
\mathbb{R}^{d}}\Gamma (r,dx)Z(r,x){\bf 1}_{\{|Z(rt,x)|\leq N\}}f(x-y)\varphi
(y)dxdydr\right)  \\
&\quad =E\left( {\bf 1}_B \int_{s}^{t}\int_{\mathbb{R}^{d}}\int_{\mathbb{R}^{d}}\Gamma
(r,dx)Z(r,x)f(x-y)\varphi (y)dxdydr\right),
\end{align*}
so we can identify $\tilde G$ with $G$. \hfill\qed

\begin{remark}
As a consequence of Proposition \ref{prop1}, we are able to define the stochastic integral of $G=Z \Gamma$ with respect to $W$:
$$G\cdot W =\int_0^T \int_{\mathbb{R}^{d}} G(s,y) W(ds,dy)=\int_0^T \int_{\mathbb{R}^{d}} \Gamma(s,y) Z(s,y) W(ds,dy).$$
Moreover, using the same ideas as in \cite{QS} (see also \cite{quer}, Theorem 1.2.5), one can obtain bounds for the $L^p(\Omega)-$norm of 
$G\cdot W$. More precisely, suppose that 
$$\sup_{(t,x)\in [0,T]\times \mathbb{R}^{d}} E(|Z(t,x)|^p)< \infty,$$
for some $p\geq 2$. Then
\begin{align}
& E(|G\cdot W|^p) = E\left( \left| \int_0^T \int_{\mathbb{R}^{d}} G(s,y) W(ds,dy)\right|^p\right) \nonumber \\
& \quad  
\leq C_p (\nu_T)^{\frac{p}{2}-1} \int_0^T \left( \sup_{x\in \mathbb{R}^{d}} E(|Z(s,x)|^p) \right) \int_{\mathbb{R}^{d}} |\mathcal{F} \Gamma (s)(\xi)|^2 \mu(d\xi) ds,
\label{pbound}
\end{align}
where
$$\nu_T =\int_0^T \int_{\mathbb{R}^{d}} |\mathcal{F} \Gamma (t)(\xi)|^2 \mu(d\xi) dt.$$
\label{Lpbound}
\end{remark}

\begin{remark}
Since the noise's correlation is of the form $(x,y)\mapsto f(x-y)$, it is natural to be interested in spatially homogeneous situations. Indeed, suppose that we add the following hypothesis on the process $Z$: for all $s\in [0,T]$ and $x,y\in \mathbb{R}^{d}$ we have
$$E(Z(s,x)Z(s,y))=E(Z(s,0)Z(s,y-x)).$$
Then, owing to \cite{Da}, p. 10, we may construct a non-negative tempered measure $\mu_s^Z$ on $\mathbb{R}^{d}$ such that
$$\|G\|_{L^2(\Omega; \mathcal{H}_T)}= E(|G\cdot W|^2)=\int_0^T \int_{\mathbb{R}^{d}} |\mathcal{F} \Gamma(s)(\xi)|^2 \mu_s^Z(d\xi) ds.$$
\label{hom}
\end{remark}

As we will see in the next section, the main examples of deterministic measures $\Gamma$ will correspond to fundamental solutions associated to second order differential operators. First, we sketch the construction of the stochastic integral in a Hilbert-valued setting, that is when the process $Z$ takes values in some Hilbert space, usually different from $\mathcal{H}$. \\

\noindent {\bf Hilbert-valued stochastic integrals}\\

Let $\mathcal{A}$ be a separable real Hilbert space with inner-product and norm denoted
by
$\langle\cdot,\cdot\rangle_{\mathcal{A}}$ and $\Vert \cdot\Vert_{\mathcal{A}}$,
respectively.
Let $K=\{K(t,x),  (t,x)\in[0,T]\times \mathbb{R}^{d}\}$
be an $\mathcal{A}-$valued predictable process satisfying the following condition:
\begin{equation}
\sup_{(t,x)\in[0,T]\times \mathbb{R}^{d}}
E\left(||K(t,x)||_{\mathcal{A}}^{2}\right)<\infty.
\label{10}
\end{equation}
Our purpose is to define the stochastic integral of elements of the form $\Gamma  K=\Gamma(t,dx)  K(t,x) \in L^2(\Omega \times  [0,T]; \mathcal{H} \otimes \mathcal{A})$.

Let $(e_{j},j\geq 0)$ be a complete orthonormal  system of $\mathcal{A}$. Set
$K^j(t,x) = \langle K(t,x),e_{j}\rangle_{\mathcal{A}}$,  $(t,x)\in [0,T]\times
\mathbb{R}^{d}$.
According to Proposition \ref{prop1}, for any $j\geq 0$ the element $G^j= G^j(t,x)= \Gamma(t,dx) K^j (t,x)$ belongs to 
$L^2(\Omega \times [0,T]; \mathcal{H})$ and, therefore, we may integrate it with respect to the noise $W$: 
$$G^j \cdot W=\int_0^T \int_{\mathbb{R}^{d}} \Gamma(s,y) K^j(s,y) W(ds,dy).$$
We define, for $G=\Gamma  K$,
$$ G\cdot W := \sum_{j\geq 0} G^j \cdot W.$$
Owing to (\ref{10}) and Proposition \ref{prop1}, it can be proved that the above series is convergent and therefore $G\cdot W$ defines an element of $L^{2}(\Omega ;\mathcal{A})$ (see also 
\cite{QS}, Remark 1). Moreover, using the same arguments as for the proof of (\ref{pbound}), we have the following bound for the moments of $G\cdot W$ in $\mathcal{A}$:
\begin{align}
& E\big(|| G\cdot W||_{\mathcal{A}}^{p}\big) \nonumber \\
& \quad \le C_p (\nu_T)^{\frac{p}{2}-1} \int_0^T  \sup_{x\in \mathbb{R}^{d}} E(||K(s,x)||^p_ {\mathcal{A}})
\int_{\mathbb{R}^{d}}  |\mathcal{F} \Gamma(s)(\xi)|^{2} \mu(d\xi) ds,
\label{11}
\end{align}
for all $p\geq 2$.


\section{Existence and uniqueness of solutions}
\label{ex}

Recall that a solution to Equation (\ref{1}) is a real-valued adapted stochastic process $u = \{u(t,x), (t,x)\in [0,T]\times \mathbb{R}^{d}\}$
satisfying
\begin{align}
u(t,x)=& \int_0^t \int_{\mathbb{R}^{d}} \Gamma(t-s,x-y) \sigma(u(s,y)) W(ds,dy) \nonumber \\
&+ \int_0^t \int_{\mathbb{R}^{d}} b(u(t-s,x-y)) \Gamma(s,dy) ds.
\label{eq}
\end{align}
We assume that   $Z=Z(s,y)=\sigma(u(s,y))$ satisfy the hypothesis of Proposition \ref{prop1}, so that the stochastic integral on the right hand-side is well-defined. 

We suppose that $\sigma$ and $b$ are real-valued Lipschitz functions. Under these conditions, we may state an existence and uniqueness of solution's theorem:

\begin{theorem}
Suppose that the fundamental solution $\Gamma$ of $Lu=0$ satisfies Hypothesis A. 
Then, Equation (\ref{eq}) has a unique solution $\{u(t,x), (t,x)\in [0,T]\times \mathbb{R}^{d}\}$ which is continuous in $L^2$ and satisfies
$$\sup_{(t,x)\in [0,T]\times \mathbb{R}^{d}} E(|u(t,x)|^p) <\infty,$$
for all $T>0$ and $p\geq 1$.
\label{sol}
\end{theorem}
For the proof we refer to \cite{Da}, Theorem 13, where the Walsh-Dalang equivalent setting is used.

Let us now enumerate some examples of differential operators whose associated fundamental solution fulfills the hypothesis of 
Theorem \ref{sol}.

\begin{example}
{\it The wave equation}. Let $\Gamma_d$ be the fundamental solution of the wave equation in $\mathbb{R}^{d}$, that is, $\Gamma_d$ is the solution of 
$$\frac{\partial^2 \Gamma_d}{\partial t^2} - \Delta \Gamma_d =0,$$
with vanishing initial conditions. It is known that for $d=1,2,3$, $\Gamma_d$ is given, respectively, by
\begin{align*}
\Gamma_1(t)&=\frac{1}{2} {\bf 1}_{\{|x|<t\}},\\
\Gamma_2 (t)& = C (t^2-|x|^2)_+^{-1/2},\\
\Gamma_3(t) &= \frac{1}{4\pi t}\sigma_t,
\end{align*}
where $\sigma_t$ denotes the surface measure on the three-dimensional sphere of radius $t$. In particular, for each $t$, $\Gamma_d(t)$ has compact support. It is important to remark that only in these cases $\Gamma_d$ defines a non-negative measure. Furthermore, for all dimensions $d\geq 1$, we have a unified expression for the Fourier transform of $\Gamma_d(t)$:
$$\mathcal{F} \Gamma_d (t)(\xi)= \frac{\sin (2\pi t |\xi|)}{2\pi |\xi|}.$$
Elementary estimates show that there are positive constants $c_1$ and $c_2$ depending on $T$ such that
$$\frac{c_1}{1+|\xi|^2}\leq \int_0^T \frac{ \sin^2 (2\pi t\xi)}{4\pi^2 |\xi|^2} dt \leq \frac{c_2}{1+|\xi|^2}.$$
Therefore, $\Gamma_d$ satisfies condition (\ref{5}) if and only if
\begin{equation}
\int_{\mathbb{R}^{d}} \frac{\mu(d\xi)}{1+|\xi|^2}<\infty.
\label{12}
\end{equation}
\label{wave}
\end{example} 

\begin{example}
{\it The heat equation}. Let $\Gamma$ be the fundamental solution of the heat equation in $\mathbb{R}^{d}$ and with vanishing initial conditions, that is
$$\frac{\partial \Gamma}{\partial t} -\frac{1}{2}\Delta \Gamma =0.$$
Then, $\Gamma$ is given by the Gaussian density:
$$\Gamma(t,x)=(2\pi t)^{-d/2} \exp\left(-\frac{|x|^2}{2 t}\right)$$
and 
$$\mathcal{F} \Gamma (t)(\xi)=\exp(-4\pi^2 t|\xi|^2).$$
Because
$$\int_0^T \exp(-4\pi^2 t|\xi|^2) dt =\frac{1}{4\pi^2 |\xi|^2} (1-\exp(-4\pi^2 T|\xi|^2)),$$
we conclude that condition (\ref{5}) holds if and only if (\ref{12}) holds.
\label{heat}
\end{example}

Let us express condition (\ref{12}) in terms of the covariance function $f$. Indeed, as it is pointed out in \cite{Da}, condition (\ref{12}) is always true when $d=1$; for $d=2$, (\ref{12}) holds if and only if
$$\int_{|x|\leq 1} f(x) \log\frac{1}{|x|} dx <\infty,$$
and for $d\geq 3$, (\ref{12}) holds if and only if
$$\int_{|x|\leq 1} f(x) \frac{1}{|x|^{d-2}} dx <\infty.$$  


\section{Existence of density}
\label{existencia}

In this section we aim to prove that the solution to Equation (\ref{3}), at any point $(t,x)\in (0,T]\times \mathbb{R}^{d}$, is a random variable whose law admits a density with respect to Lebesgue measure on $\mathbb{R}$. For this, we will make use of the techniques provided by the Malliavin calculus and, more precisely, we will apply Bouleau-Hirsch's criterion (see \cite{bh} or Theorem 2.1.2 in \cite{nualart}).\\

First of all, we describe the Gaussian context in which we will use the tools of the Malliavin calculus. Namely, we consider the Hilbert space $\mathcal{H}_T=L^2([0,T];\mathcal{H})$ and the Gaussian family of random variables 
$(W(h), h\in \mathcal{H}_T)$ defined at the very beginning of Section \ref{integral}.
Then
$(W(h), h\in \mathcal{H}_T)$ is a centered Gaussian process such that $E(W(h_1)W(h_2))=\langle h_1,h_2\rangle_{\mathcal{H}_T}$, 
$h_1,h_2\in \mathcal{H}_T$,
and we can use the differential Malliavin calculus based on it (see, for instance, \cite{nualart}). The Malliavin derivative is denoted by $D$ and, for any $N\geq 1$, the domain of the iterated derivative $D^N$ in $L^p(\Omega; \mathcal{H}_T^{\otimes N})$ is denoted by 
$\mathbb{D}^{N,p}$, for any $p\geq 2$. We shall also use the notation 
$$\mathbb{D}^\infty =\cap_{p\geq 1}\cap_{k\geq 1} \mathbb{D}^{k,p}.$$

The first step in order to apply Bouleau-Hirsch's criterion is to study the Malliavin differentiability of $u(t,x)$, for all fixed $(t,x)\in (0,T]\times \mathbb{R}^{d}$. Recall that, for any random variable $X$ in the domain of the derivative operator $D$, $DX$ defines an $\mathcal{H}_T-$valued random variable. In particular, for some fixed $r\in [0,T]$, $DX(r)$ is an element of $\mathcal{H}$, which will be denoted by $D_rX$. In the sequel we will use the notation $\cdot$ and $*$ to denote, respectively, the time and $\mathcal{H}$ variables.

\begin{proposition}
Assume that $\Gamma$ satisfies Hypothesis A. Suppose also that the coefficients $b$ and $\sigma$ are $\mathcal{C}^1$ functions with bounded Lipschitz continuous derivatives. Then, for any $(t,x)\in [0,T]\times \mathbb{R}^{d}$, $u(t,x)$ belongs to 
$\mathbb{D}^{1,p}$, for any $p\in [1,\infty)$. 

Moreover, the Malliavin derivative $Du(t,x)$ defines an $\mathcal{H}_T-$valued process that satisfies the following linear stochastic differential equation: 
\begin{align}
D_r u(t,x) = & \sigma(u(r,*))\Gamma(t-r,x-*) \nonumber\\
& + \int_r^t \int_{\mathbb{R}^{d}} \Gamma(t-s,x-y) \sigma'(u(s,y)) D_ru(s,y)W(ds,dy)\nonumber\\
& + \int_r^t \int_{\mathbb{R}^{d}}  b'(u(s,x-y)) D_r u(s,x-y)\Gamma(t-s,dy) ds,
\label{eqmal}
\end{align}
for all $r\in [0,T]$.
\label{difmal}
\end{proposition}

The stochastic integral on the right hand-side of Equation (\ref{eqmal}) must be understood by means of the Hilbert-valued integration setting described at the very final part of Section \ref{integral}.

Concerning the Hilbert-valued pathwise integral, it is defined as follows: let $\mathcal{A}$ be a Hilbert space, $(e_j)_{j\geq 1}$ a complete orthonormal system of $\mathcal{A}$ and $\{Y(s,y), (s,y)\in [0,T]\times \mathbb{R}^{d}\}$ an $\mathcal{A}-$valued stochastic process such that 
$$\sup_{(s,y)\in [0,T]\times \mathbb{R}^{d}} E(\|Y(s,y)\|^2_{\mathcal{A}})<+\infty.$$ 
Then, the $\mathcal{A}-$valued integral 
$$\mathcal{I}_t=\int_0^t  \int_{\mathbb{R}^{d}} Y(s,y)\Gamma(s,dy) ds$$
is determined by the components $\left( \int_0^t  \int_{\mathbb{R}^{d}} \langle Y(s,y),e_j\rangle_{\mathcal{A}}\; \Gamma(s,dy) ds, j\geq 1\right)$, which are real-valued integrals. Moreover, one can obtain an upper bound for the moments of the above integral (see \cite{quer}, p. 24):
\begin{equation}
E(|\mathcal{I}_t|^p)\leq \int_0^t \sup_{z\in \mathbb{R}^{d}} E(\|Y(s,z)\|_\mathcal{A}^p \int_{\mathbb{R}^{d}} \Gamma(s,dy) ds,\; p\geq 2.
\label{deter}
\end{equation}
Eventually, notice that owing to Proposition \ref{prop1}, the first term on the right hand-side of (\ref{eqmal}) is fully defined.

\medskip
\noindent {\it Proof of Proposition \ref{difmal}}. The statement is almost an immediate consequence of Theorem 2 in \cite{QS}.  Indeed, the authors of this latter reference prove that, under the standing hypothesis and for any fixed $(t,x)\in [0,T]\times \mathbb{R}^{d}$, the random variable $u(t,x)$     
belongs to $\mathbb{D}^{1,p}$, for any $p\in [1,\infty)$. In addition, they show that there exists an $\mathcal{H}_T-$valued stochastic process 
$\{\Theta(t,x), (t,x)\in [0,T]\times \mathbb{R}^{d}\}$ satisfying 
$$\sup_{(t,x)\in [0,T]\times \mathbb{R}^{d}} E(\|\Theta(t,x)\|_{\mathcal{H}_T}^p)<\infty$$
and such that, in $\mathcal{H}_T$,
\begin{align*}
D u(t,x) = &
\Theta(t,x)+ \int_0^t \int_{\mathbb{R}^{d}} \Gamma(t-s,x-y) \sigma'(u(s,y)) Du(s,y)W(ds,dy)\\
& + \int_0^t \int_{\mathbb{R}^{d}}  b'(u(s,x-y)) D u(s,x-y)\Gamma(t-s,dy) ds. 
\end{align*}
Moreover, it holds that
$$E(\|\Theta(t,x)\|_{\mathcal{H}_T}^2)=E(\| \Gamma(t-\cdot,x-*) \sigma(u(\cdot,*))\|_{\mathcal{H}_T}^2).$$
Hence, in order to conclude the proof, we only need to  show that the Hilbert-valued random variables $\Theta(t,x)$ and 
$\Gamma(t-\cdot,x-*) \sigma(u(\cdot,*))$ coincide  as elements of $\mathcal{H}_T$. 

Let $(\Gamma_n)_{n\geq 1}$ be the family of smooth functions defined in the proof of Lemma \ref{ggamma}. Then, 
in the proof of Theorem 2 in \cite{QS} the process $\Theta$ is defined by the following limit in $\mathcal{H}_T$:
$$\Theta(t,x)=\mathcal{H}_T-\lim_{n\rightarrow \infty} \Gamma_n(t-\cdot,x-*) \sigma(u_n(\cdot,*)),$$
where $\{u_n(t,x), (t,x)\in [0,T]\times \mathbb{R}^{d}\}$ is the unique mild solution to an equation of the form (\ref{3}) but replacing $\Gamma$ by $\Gamma_n$. 

As a consequence of the proof of Proposition 3 from \cite{QS}, it is readily checked that
$$
\lim_{n\rightarrow 0}  E(\| \Gamma_n(t-\cdot,x-*)[\sigma(u_n(\cdot,*))-\sigma(u(\cdot,*))]\|_{\mathcal{H}_T}^2)=0
$$
and 
$$\lim_{n\rightarrow 0}  E(\| [\Gamma_n(t-\cdot,x-*)-\Gamma(t-\cdot,x-*)]\sigma(u(\cdot,*))\|_{\mathcal{H}_T}^2)=0.$$
Thus, we get that $\Theta(t,x)=\sigma(u(t-\cdot,x-*)) \Gamma(t-\cdot, x-*)$. \hfill \qed

\vspace{0.3cm}

The main result of the section is the following:

\begin{theorem}
Assume that $\Gamma$ satisfies Hypothesis A. Suppose also that the coefficients $b$ and $\sigma$ are $\mathcal{C}^1$ functions with bounded Lipschitz continuous derivatives and that $|\sigma(z)|\geq c>0$, for all $z\in \mathbb{R}$ and some positive constant $c$. Then, for all $t>0$ and $x\in \mathbb{R}^{d}$, the random variable $u(t,x)$ has an absolutely continuous law with respect to Lebesgue measure on $\mathbb{R}$.
\label{densitat}
\end{theorem}

\noindent
{\it Proof}. Owing to Bouleau-Hirsch's criterion and Proposition \ref{difmal}, it suffices to show that $\|Du(t,x)\|_{\mathcal{H}_T}>0$ almost surely.

To begin with,  from Equation (\ref{eqmal})  we obtain
\begin{equation}
\int_0^t \|D_s u(t,x)\|_\mathcal{H}^2 ds\geq \frac{1}{2}\int_{t-\delta}^t \|\Gamma(t-s,x-*)\sigma(u(s,*))\|_\mathcal{H}^2 ds - I(t,x;\delta),
\label{15}
\end{equation}
for any $\delta>0$ sufficiently small, where 
\begin{align}
I(t,x;\delta)= & \int_{t-\delta}^t \left\| \int_s^t \int_{\mathbb{R}^{d}} \Gamma(t-r,x-z) \sigma'(u(r,z)) D_s u(r,z) W(dr,dz) \right.\nonumber\\
&\quad + \left. \int_s^t \int_{\mathbb{R}^{d}}\Gamma(t-r,dz) b'(u(r,x-z)) D_s u(r,x-z) dr\right\|_\mathcal{H}^2 ds.
\label{15.5}
\end{align} 
The above term $I(t,x;\delta)$ may be bounded by $2(I_1(t,x;\delta)+I_2(t,x;\delta)$, with
\begin{align}
I_1(t,x;\delta)= &  \int_0^\delta \left\| \int_{t-s}^t \int_{\mathbb{R}^{d}} \Gamma(t-r,x-z) \sigma'(u(r,z)) D_{t-s} u(r,z) W(dr,dz) \right\|_\mathcal{H}^2 ds, \label{15.6}\\
I_2(t,x;\delta)= & \int_0^\delta \left\| \int_{t-s}^t \int_{\mathbb{R}^{d}}\Gamma(t-r,dz) b'(u(r,x-z)) D_{t-s} u(r,x-z) dr\right\|_\mathcal{H}^2 ds.
\label{15.7}
\end{align} 
In order to bound from below the term in the left hand-side of (\ref{15}), let us first
obtain a lower bound for the first one on the right hand-side. For this, we will make use of the family of smooth functions $(\Gamma_n)_n$ and 
$(J_n^{t,x})_n$, considered in the proof of Lemma \ref{ggamma}, that
regularise the measures $\Gamma$ and $\Gamma(\cdot,x-*)\sigma(u(t-\cdot,*))$, respectively. Then, by the proof of Lemma \ref{ggamma}, the very definition of the norm in $\mathcal{H}_\delta$ and the non-degeneracy assumption on $\sigma$, we have
\begin{align}
& \int_{t-\delta}^t \|\Gamma(t-s,x-*)\sigma(u(s,*))\|_\mathcal{H}^2 ds = \|\Gamma(\cdot,x-*)\sigma(u(t-\cdot,*))\|_{\mathcal{H}_\delta}^2 \nonumber\\
&\quad = \lim_{n\rightarrow \infty} \|  J_n^{t,x} \|_{\mathcal{H}_\delta}^2\nonumber\\
& \quad = \lim_{n\rightarrow \infty} \int_0^\delta \int_{\mathbb{R}^{d}}\int_{\mathbb{R}^{d}} J_n^{t,x}(s,y) f(y-z) J_n^{t,x}(s,z)
dydz ds\nonumber \\
&\quad \geq c^2 \lim_{n\rightarrow \infty} \int_0^\delta \int_{\mathbb{R}^{d}}\int_{\mathbb{R}^{d}} \Gamma_n(s,x-y) f(y-z)\Gamma_n(s,x-z)\nonumber \\
& \quad = c^2 \lim_{n\rightarrow \infty} \|\Gamma_n(\cdot,x-*) \|_{\mathcal{H}_\delta}^2 = c^2 \|\Gamma(\cdot,x-*)\|_{\mathcal{H}_\delta}^2=c^2 g(\delta),
\label{16}
\end{align}
where 
\begin{equation} \label{gedelta}
g(\delta):=\int_0^\delta \int_{\mathbb{R}^{d}} |\mathcal{F} \Gamma(s)(\xi)|^2 \mu(d\xi) ds. 
\end{equation}

Now we find out upper bounds for the expectation of the terms $I_1(t,x;\delta)$ and $I_2(t,x;\delta)$. First, in order to deal with 
the former term, one can use the bound (\ref{11}). Thus, taking into account that $\sigma'$ is bounded, we get the following estimate:
\begin{align}
& E(I_1(t,x;\delta))\nonumber \\ 
&\quad  \leq C  \sup_{(\tau,y)\in (0,\delta)\times \mathbb{R}^{d}} E\left( \|D_{t-\cdot} u(t-\tau,y)\|_{\mathcal{H}_\delta}^2\right) \int_0^\delta \int_{\mathbb{R}^{d}} |\mathcal{F} \Gamma(s)(\xi)|^2 \mu(d\xi) ds\nonumber \\
& \quad= C \sup_{(\tau,y)\in (0,\delta)\times \mathbb{R}^{d}} E\left( \|D_{t-\cdot} u(t-\tau,y)\|_{\mathcal{H}_\delta}^2\right) g(\delta).
\label{17}
\end{align}
On the other hand, by (\ref{deter}) the term $E(I_2(t,x;\delta))$ corresponding to the Hilbert-valued pathwise integral can be bounded by 
\begin{equation} 
E(I_2(t,x;\delta))\leq C \sup_{(\tau,y)\in (0,\delta)\times \mathbb{R}^{d}} E\left( \|D_{t-\cdot} u(t-\tau,y)\|_{\mathcal{H}_\delta}^2\right) h(\delta),
\label{17.2}
\end{equation}
where $h(\delta)=\int_0^\delta \int_{\mathbb{R}^{d}} \Gamma(s,dy) ds$.\\

At this point, we will make use of the following fact (see Lemma 5 in \cite{qs2} and \cite{quer}, p. 53):  
\begin{equation}
\sup_{(\tau,y)\in (0,\delta)\times \mathbb{R}^{d}} E\left( \|D_{t-\cdot} u(t-\tau,y)\|_{\mathcal{H}_\delta}^{2q}\right)
\leq C (g(\delta))^q,
\label{17.5}
\end{equation}
for any $q\geq 1$. Hence, by (\ref{17}) and (\ref{17.2}) the terms $E(I_1(t,x;\delta))$ and $E(I_2(t,x;\delta))$ may be bounded, up to constants, respectively, by $(g(\delta))^2$ and $g(\delta) h(\delta)$, which implies that
\begin{equation}
E(I(t,x;\delta))\leq C g(\delta)(g(\delta)+h(\delta)).
\label{18}
\end{equation}
For any fixed small $\delta>0$, let $n$ be a sufficiently large positive integer such that $\frac{1}{n}\leq \frac{c^2}{2}g(\delta)$. Then, owing to 
(\ref{16}) and (\ref{18}) and applying Chebyshev's inequality, we obtain
\begin{align}
P\left( \int_0^t \|D_s u(t,x)\|_\mathcal{H}^2 ds <\frac{1}{n}\right) & \leq P\left( I(t,x;\delta)\geq \frac{c^2}{2} g(\delta)-\frac{1}{n}\right) \nonumber\\
& \leq \left(\frac{c^2}{2}g(\delta)-\frac{1}{n}\right)^{-1} E(I(t,x;\delta))\nonumber \\
& \leq \left(\frac{c^2}{2}g(\delta)-\frac{1}{n}\right)^{-1} g(\delta)(g(\delta)+h(\delta)).
\label{19}
\end{align}
Therefore
$$\lim_{n\rightarrow \infty} P\left( \int_0^t \|D_s u(t,x)\|_\mathcal{H}^2 ds <\frac{1}{n}\right) \leq C (g(\delta)+h(\delta)),$$
and the latter term converges to zero as $\delta$ tends to zero. Hence,
$$P\left( \int_0^t \|D_s u(t,x)\|_\mathcal{H}^2 ds =0\right)=0,$$
which concludes the proof.   \hfill\qed

\begin{remark}
In the particular case of the three-dimensional stochastic wave equation, 
the above Theorem \ref{densitat} generalises  Theorem 3 in the reference \cite{QS}. 
\end{remark}


\section{Smoothness of the density}
\label{regularitat}

This section is devoted to prove that, for any fixed $(t,x)\in (0,T]\times \mathbb{R}^{d}$, the law of the random variable $u(t,x)$ has an infinitely differentiable density with respect to Lebesgue measure on $\mathbb{R}$. This will be achieved by showing that $u(t,x)$ belongs to the space $\mathbb{D}^\infty$ and that the inverse of the Malliavin matrix of $u(t,x)$ has moments of all order (see, for instance, Theorem 2.1.4 in \cite{nualart}).   

Recall that for any {\it differentiable} random variable $X$ and any $N\geq 1$, the iterated Malliavin derivative $D^N X$ defines an element of the Hilbert space $L^2(\Omega; \mathcal{H}_T^{\otimes N})$. As for the case $N=1$, for any $r=(r_1,\dots,r_N)\in [0,T]^N$, the element $D X(r)$ of $\mathcal{H}^{\otimes N}$ will be denoted by $D_rX$. We will also use the notation 
$$D^{N}_{((r_{1},\varphi_{1}),\dots,(r_{N},\varphi_{N}))}X=
\langle D^{N}_{(r_{1},\dots,r_{N})} X, \varphi_{1}\otimes \dots \otimes
\varphi_{N}
\rangle_{\mathcal{H}^{\otimes N}},$$
for $r_{i}\in [0,T]$, $\varphi_{i}\in\mathcal{H}$, $i=1,\dots,N$.
In particular, we have that
$$
\|D^{N} X\|^{2}_{\mathcal{H}_T^{\otimes N}}=\int_{[0,T]^{N}}dr_{1}\dots dr_{N}
\sum_{j_{1},\dots,j_{N}}
|D_{((r_{1},e_{j_{1}}),\dots,(r_{N},e_{j_{N}}))} X|^{2}, 
$$
where $(e_{j})_{j\geq 0}$ is a complete orthonormal system of $\mathcal{H}$. Let
$$\Delta_{\alpha}^N (g,X):= D^{N}_{\alpha} g(X) - g'(X) D_{\alpha}^{N} X,$$
where $\alpha=((r_1,\varphi_1),\dots,(r_N,\varphi_N))$, $r_i\in [0,T]$ and $\varphi_i\in \mathcal{H}$.
Notice that $\Delta_{\alpha}^N (g,X)=0$ if $N=1$ and it only depends on the
Malliavin 
derivatives up
to the order $N-1$ if $N>1$.\\

We now state the main result concerning the Malliavin regularity of the solution $u(t,x)$.

\begin{proposition}
Assume that $\Gamma$ satisfies Hypothesis A. Suppose also that the coefficients $\sigma$ and $b$ are $\mathcal{C}^\infty$ functions with bounded derivatives of any order greater than or equal to one. Then, for every $(t,x)\in [0,T]\times \mathbb{R}^{d}$, the random variable $u(t,x)$ belongs to the space $\mathbb{D}^\infty$. 

The iterated Malliavin derivative $D^Nu(t,x)$ satisfies the following equation in $L^p(\Omega; \mathcal{H}_T^{\otimes N})$, for any $p\geq 1$ and 
$N\geq 1$:
\begin{align}
& D^{N}u(t,x)  = Z^{N}(t,x) \nonumber \\
& + \int_{0}^{t}\int_{\mathbb{R}^{d}} \Gamma(t-s,x-z) [\delta^N(\sigma,u(s,z)) 
+ D^{N}u(s,z)\sigma'(u(s,z))] W(ds,dz) \nonumber \\
& + \int_{0}^{t}ds \int_{\mathbb{R}^{d}} \Gamma(s,dz) [\delta^N(b,u(t-s,x-z))\nonumber \\
&\quad \quad 
+ D^{N}u(t-s,x-z)b'(u(t-s,x-z))], 
\label{21}
\end{align}
where $Z^N(t,x)$ is the element of $L^p(\Omega; \mathcal{H}_T^{\otimes N})$ defined by
\begin{equation}
\langle Z^N_r(t,x), e_{j_1}\otimes\dots \otimes e_{j_N}\rangle_{\mathcal{H}^{\otimes N}} =
\sum_{i=1}^{N} \langle \Gamma(t-r_i,x-*)
D^{N-1}_{\hat{\alpha}_i}
\sigma(u(r_i,*)),e_{j_i}\rangle_{\mathcal{H}},
\label{ci}
\end{equation}
for any $r=(r_1,\dots,r_N)\in [0,T]^N$ and $j_1,\dots,j_N\in \{1,\dots,N\}$.

Moreover, it holds that
$$\sup_{(s,y)\in [0,T]\times \mathbb{R}^{d}} E(\|D^{N}u(s,y)\|^{p}_{\mathcal{H}_T^{\otimes N}}
)<+\infty,$$
for all $p\ge 1$.
\label{malinf}
\end{proposition}

\noindent
{\it Proof}. It is almost an immediate consequence of Theorem 1 in \cite{qs2} and Proposition \ref{difmal} from the preceding Section 
\ref{existencia}. 

Namely, we just need to check that, using the same notation as in the proof of Theorem 1 in \cite{qs2}, the sequence of 
$L^2(\Omega; \mathcal{H}_T^{\otimes N})-$valued random variables $(Z^{N,n}(t,x))_{n\geq 1}$ converges to $Z^N(t,x)$ as $n$ tends to infinity. We should mention that in that reference, $Z^{N,n}(t,x)$ is constructed using a regularisation procedure, that is smoothing the measure $\Gamma$, and by means of a similar expression to (\ref{ci}). In \cite{qs2}, the objective was to define the initial condition of the linear stochastic equation satisfied by the iterated Malliavin derivative  $D^N u(t,x)$ as the limit of $Z^{N,n}(t,x)$. We claim that this limit equals to $Z^N(t,x)$, defined in the present proposition's statement (see (\ref{ci})).

Indeed, the convergence of $Z^{N,n}(t,x)$ to $Z^N(t,x)$ in $L^2(\Omega; \mathcal{H}_T^{\otimes N})$ can be easily studied using the same arguments as in the proof of Lemma 3 in \cite{qs2}.   \hfill \qed

\vspace{0.3cm}

We are now in position to state and prove the main result of the paper. 

\begin{theorem}
Assume that $\Gamma$ satisfies Hypothesis A, the coefficients $\sigma$ and $b$ are $\mathcal{C}^\infty$ functions with bounded derivatives of any order greater than or equal to one and $|\sigma(z)|\geq c>0$, for all $z\in \mathbb{R}$.
Moreover, suppose that there exist $\gamma >0$ such that for all $\tau\in (0,1]$,
\begin{equation}
   \int_0^\tau \int_{\mathbb{R}^{d}} |\mathcal{F} \Gamma(s)(\xi)|^2 \mu(d\xi) ds
   \ge C_1 \tau^{\gamma},
\label{lowbound}
\end{equation} 
for some positive constant  $C_1$. Then, for all $(t,x)\in [0,T]\times \mathbb{R}^{d}$, the law of $u(t,x)$ has a $\mathcal{C}^\infty$ density with respect to Lebesgue measure on $\mathbb{R}$.
\label{regularitatdensitat}
\end{theorem}

\noindent
{\it Proof.} In view of Proposition \ref{malinf}, we need to show that the inverse of the Malliavin matrix of $u(t,x)$ has moments of all order, that is
$$E\left( \left| \int_0^T \|D_s u(t,x)\|_\mathcal{H}^2  ds \right|^{-q}  \right)<+\infty,$$
for all $q\geq 2$. 
It turns out (see, for instance, Lemma 2.3.1 in \cite{nualart}) that it suffices to check that for any $q\geq 2$,
there exists an $\varepsilon _{0}(q)>0$ such that  for all $\varepsilon \leq
\varepsilon _{0}$
\begin{equation}
P\left( \int_{0}^{t}\left\| D_s u(t,x) \right\|_{\mathcal{H}}^{2} ds<\varepsilon \right) \leq C \varepsilon ^{q}.
\end{equation}
Proceeding as in the proof of Theorem \ref{densitat}, for any $\delta>0$ sufficiently small we obtain the following estimate:
\begin{align}
P\left( \int_0^t \|D_s u(t,x)\|^2_\mathcal{H} ds<\varepsilon\right) & \leq  P\left( I(t,x;\delta) \geq \frac{c^2}{2} g(\delta)-\varepsilon\right)\nonumber\\
& \leq \left( \frac{c^2}{2}g(\delta) -\varepsilon\right)^{-p} E(|I(t,x;\delta)|^p),
\label{23}
\end{align}
for any $p>0$, where we recall that $I(t,x;\delta)$ is defined by (\ref{15.5}) and 
$g(\delta)$ is given by (\ref{gedelta}).  

We decompose now the term $I(t,x;\delta)$ as in the proof of Theorem \ref{densitat}, so that we need to find upper bounds for $E(|I_i(t,x;\delta)|)^p$, $i=1,2$ (see (\ref{15.6}) and (\ref{15.7})). On one hand, owing to H\"older's inequality and (\ref{11}) we get
\begin{align}
& E(|I_1(t,x;\delta)|^p)\nonumber  \\
&\quad  = E\left( \int_0^\delta \left\| \int_{t-s}^t \int_\mathbb{R}^{d} \Gamma(t-r,x-z)\sigma'(u(r,z)) D_{t-s} u(r,z)W(dr,dz) \right\|_\mathcal{H}^2 ds\right)^p \nonumber \\
& \quad \leq \delta^{p-1} E\left( \int_0^\delta \left\| \int_{t-s}^t \int_\mathbb{R}^{d} \Gamma(t-r,x-z)\sigma'(u(r,z)) D_{t-s} u(r,z)W(dr,dz) \right\|_\mathcal{H}^{2p} ds\right) \nonumber\\
&\quad \leq \delta^{p-1} (g(\delta))^p \sup_{(\tau,y)\in [0,T]\times \mathbb{R}^{d}} E\left( \|D_{t-\cdot} u(t-\tau,y)\|_{\mathcal{H}_T}^{2p}\right).
\label{23.4}
\end{align}
The above estimate (\ref{23.4})   let us conclude that
$$
E(|I_1(t,x;\delta)|^p) \leq C \delta^{p-1} (g(\delta))^{p}.
$$
On the other hand, using similar arguments but for the  Hilbert-valued pathwise integral (see (\ref{deter})), one proves that $E(|I_2(t,x;\delta)|^p)$ may be bounded, up to some positive constant, by $\delta^{p-1} (g(\delta))^p  $.
Thus, we have proved that
\begin{equation}
 P\left( \int_0^t \|D_s u(t,x)\|^2_\mathcal{H} ds  < \epsilon\right) 
  \leq C \left( \frac{c^2}{2}g(\delta) -\epsilon\right)^{-p} \delta^{p-1} (g(\delta))^p .
\label{24}
\end{equation}
At this point, we choose $\delta=\delta(\epsilon)$ in such a way that $g(\delta)=\frac{4}{c^2} \epsilon$. By (\ref{lowbound}), this implies that 
$\frac{4}{c^{2}}\varepsilon \geq C\delta
^{\gamma }$, that is $\delta \leq C\varepsilon ^{\frac{1}{\gamma }}$. Hence,%
\[
P\left( \int_{0}^{t}\left\| D_{s}u(t,x)\right\| _{\mathcal{H}%
}^{2}ds<\varepsilon \right) \leq C\varepsilon ^{\frac{p-1}{\gamma }},
\]%
and it suffices to take $p$ sufficiently large such that $  \frac{p-1}{\gamma } \geq q$.
\hfill \qed

\begin{remark}
As it is pointed out in \cite{QS}, Appendix A (see also \cite{lev}), when $\Gamma$ is the fundamental solution of the wave equation in $\mathbb{R}^d$, with 
$d=1,2,3$, then condition  (\ref{lowbound}) is satisfied with $\gamma=3$.  

On the other hand, if $\Gamma$ is the fundamental solution of the heat equation on $\mathbb{R}^d$, $d\geq 1$, then 
condition (\ref{lowbound}) is satisfied  with  any $\gamma\geq 1$   (see Lemma 3.1 in \cite{mms}).  
\end{remark}

\begin{remark}
Theorem \ref{regularitatdensitat} provides a generalisation of Theorem 3 in \cite{qs2} for the case of the three-dimensional wave equation (see also \cite{marta}). Moreover, it also generalises the results in \cite{mms} for the stochastic wave equation with space dimension $d=1,2$ and the stochastic heat equation in any space dimension.   
\end{remark}



\begin{thebibliography}{99}

\bibitem{bh} Bouleau, N. and Hirsch, F., Dirichlet Forms and
Analysis on Wiener space, de Gruyter Studies in Mathematics {\bf 14},
Walter de Gruyter, Berlin, New York 1991.

\bibitem{carmona} Carmona, R. and Nualart, D.,
{\it{Random nonlinear wave equations: Smoothness of the solutions}},
Probab. Theory Relat. Fields {\bf{79}}, No.4, 469-508 (1988).

\bibitem{Da}Dalang, R. C., {\it{Extending martingale measure stochastic integral with applications to spatially homogeneous S. P. D. E's}},
Electron. J. Probab. {\bf{4}}, Paper No.6, 29 p. (1999) (electronic).

\bibitem{DF} Dalang, R. C. and Frangos, N. E., 
{\it{The stochastic wave equation in two spatial dimensions}},  
Ann. Probab.  {\bf{26}},  no. 1, 187-212   (1998).

\bibitem{dz} Da Prato, G. and Zabczyk, J., Stochastic equations in infinite dimensions. Encyclopedia of Mathematics and its Applications, 44. Cambridge University Press, Cambridge, 1992.

\bibitem{lev} O. L\'ev\^eque, Hyperbolic Stochastic Partial
Differential Equations Driven by
Boundary noises, Th\`{e}se {\bf{2452}} (2001) EPFL, Lausanne.

\bibitem{mms} M\'arquez-Carreras, D., Mellouk, M. and Sarr\`a, M.,
{\it On stochastic partial differential equations with spatially correlated
noise: smoothness of the law}, Stoch. Proc. Appl. {\bf 93}, 269-284 (2001).


\bibitem{MS} Millet, A. and Sanz-Sol\'e, M., 
{\it{A stochastic wave equation in two space dimensions: smoothness of the law}},
Ann. Probab. {\bf{27}}, No.2, 803-844 (1999). 

\bibitem{nualart} Nualart, D., The Malliavin Calculus and Related Topics, Second edition. Probability and its Applications (New York). Springer-Verlag, Berlin, 2006.

\bibitem{Pe} Peszat, S., {\it{The Cauchy problem for a nonlinear stochastic wave equation in any dimension}},  J. Evol. Equ.  {\bf{2}}, no. 3, 383-394 (2002).

\bibitem{PZ} Peszat, S. and  Zabczyk, J., {\it{Nonlinear stochastic wave and heat equations}},  Probab. Theory Related Fields  {\bf{116}},  no. 3, 421-443 (2000).

\bibitem{quer} Quer-Sardanyons, L., The stochastic wave equation: study of the law and approximations, PhD-Thesis Universitat de Barcelona 2005.

\bibitem{QS} Quer-Sardanyons, L. and Sanz-Sol\'e, M., {\it{Absolute continuity of the law of the solution to the 3-dimensional stochastic wave equation}},
J. Funct. Anal. {\bf{206}}, No.1, 1-32 (2004).

\bibitem{qs2} Quer-Sardanyons, L. and Sanz-Sol\'e, M., {\it{A stochastic wave equation in dimension 3: Smoothness of the law}},
Bernoulli {\bf{10}}, No.1, 165-186 (2004). 

\bibitem{marta} Sanz-Sol\'e, M., 
Malliavin calculus. With applications to stochastic partial differential equations. Fundamental Sciences. EPFL Press, Lausanne; distributed by CRC Press, Boca Raton, FL, 2005.

\bibitem{walsh} Walsh, J. B., {\it{An introduction to stochastic partial differential equations}},
\'Ecole d'\'et\'e de probabilit\'es de Saint-Flour XIV - 1984, Lect. Notes Math. {\bf{1180}}, 265-437 (1986).


\end{thebibliography}
\end{document}